\documentstyle[11pt]{article}
\begin{document}

\font\BBfont=msbm10
\def\Bbb#1{\hbox {\BBfont #1}} 

\newtheorem{defi}{Definition}[section]

\newcommand{\qed}{\hfill$\square$} 
\newcommand{\pr}{{\bf Proof. }} 
\newcommand{\rr}{\Bbb R}
\newcommand{\ccl}{{\cal L}}
\newcommand{\ccb}{{\cal B}}
\newcommand{\cca}{{\cal A}}
\newcommand{\ccf}{{\cal F}}
\newcommand{\ccr}{{\cal R}}
\newcommand{\nn}{\Bbb N}
\newcommand{\la}{\lambda}
\newcommand{\lz}{\left|\!\left|}
\newcommand{\rz}{\right|\!\right|}
\newcommand{\ra}{\rightarrow}
\newcommand{\vep}{\varepsilon}
\newcommand{\beq}{\begin{equation}}
\newcommand{\eeq}{\end{equation}}
\newcommand{\beqn}{\begin{eqnarray}}
\newcommand{\eeqn}{\end{eqnarray}}
\newcommand{\beqnn}{\begin{eqnarray*}}
\newcommand{\eeqnn}{\end{eqnarray*}}
\newcommand{\smallsetminus}{\setminus}
\newcommand{\varnothing}{\emptyset}

\title{\bf Ramsey dichotomies with ordinal index}
         
\author{Vassiliki Farmaki}

\date{}
\maketitle

\begin{abstract}
A system of uniform families on an infinite subset $M$ of $\nn$ is a
collection $(\cca_{\xi})_{\xi<\omega_1}$ of families of finite subsets of
$\nn$ (where, $\cca_k$ consists of all $k$--element subset of $M$, for
$k\in \nn$) with the properties that each $\cca_{\xi}$ is thin (i.e. it
does not contain proper initial segments of any of its element) and the
Cantor--Bendixson index, defined for $\cca_{\xi}$, is equal to $\xi+1$ and
stable when we restrict ourselves to any subset of $M$. We indicate how to
extend the generalized Schreier families to a system of uniform families.

Using that notion we establish the correct (countable) ordinal index
generalization of the classical Ramsey theorem (which corresponds to the
finite ordinal indices). Indeed, for a family $\ccf$ of finite subsets of
$\nn$, we obtain the following:
\begin{enumerate}
\item [(i)] For every infinite subset $M$ of $\nn$ and every countable
ordinal $\xi$, there is an infinite subset $L$ of $M$ such that either
$\cca_{\xi}\cap [L]^{<\omega}\subseteq\ccf$ or $\cca_{\xi}\cap
[L]^{<\omega}\subseteq [\nn]^{<\omega}\smallsetminus\ccf$;\\
(where $[L]^{<\omega}$ denotes the family of all finite subsets of $L$).
\item [(ii)] If, in addition $\ccf$ is hereditary and pointwise closed,
then for every infinite subset $M$ of $\nn$ there is a countable ordinal
number $\xi$ such that:
\begin{enumerate}
\item [(a)] For every ordinal number $\zeta$ with $\zeta+1<\xi$ there is an
infinite subset $L$ of $M$ such that $\cca_{\xi}\cap
[L]^{<\omega}\subseteq\ccf$. 
\item [(b)] For every ordinal number $\zeta$ with $\xi<\zeta+1$ there is an
infinite subset $L$ of $M$ such that $\ccf\cap [L]^{<\omega}\subseteq
(\cca_{\zeta})^{*}\smallsetminus \cca_{\zeta}$; which gives 
$\cca_{\xi}\cap
[L]^{<\omega}\subseteq [\nn]^{<\omega}\smallsetminus\ccf$;\\
(where generally $\cca^{*}$ denotes the family of all initial segments of
elements of $\cca$).
\item [(c)] For $\zeta=\xi+1$, both alternatives ((a) and (b)) may
materialize. 
\end{enumerate}
\item [(iii)] If $\ccf$ is hereditary, then $\ccf$ is not closed if and
only if there is an infinite subset $M$ of $\nn$ such that
$[M]^{<\omega}\subseteq \ccf$. 
\end{enumerate}
\end{abstract}


\newpage
\setcounter{section}{-1}
\section{Introduction}

Our aim, in the present paper, is to establish the proper context, and the
correct (countable) ordinal (Cantor--Bendixson type) index generalization
of the classical Ramsey theorem \cite{14} (stating that for every family
$\ccf$ of finite subsets of $\nn$, every natural number $k$ and every
infinite subset $M$ of $\nn$, there is an infinite subset $L$ of $M$, such
that all subsets of $M$ consisting of exactly $k$--elements are either in
$\ccf$ or in the complement of $\ccf$).

In this ordinal index context, the index of the classical Ramsey theorem is
a natural number, while the infinitary Galvin--Prikry theorem, or infinite
Ramsey, as is sometimes loosely called, (\cite{11}, \cite{8}, \cite{16} and
\cite{5b}) corresponds to the limiting $\omega_1$--ordinal index.

Using the notion of a uniform family given by Pudl\'{a}k and R\"{o}dl in
\cite{13} we introduce the notion of a system of uniform families
(Definition 1.3). A system of uniform families on $M$ (for $M\in [\nn]$) is
a collection $(\cca_{\xi})_{\xi<\omega_1}$ of families of finite subsets of
$M$ (with $\cca_k=[M]^k$ for $k\in
\nn$) with the properties: (i) each $\cca_{\xi}$ is thin (i.e. it does not
contain proper initial segments of any of its elements) and (ii) the
Cantor--Bendixson index defined for $\cca_{\xi}$ is precisely equal to
$\xi+1$ and does not decrease, but on the contrary is stable, when we
restrict ourselves to any infinite subset of $M$.

Every system of uniform families on $M$ is characterized by the choice, for
each countable limit ordinal number $\xi$, of an increasing sequence
$(\xi_m)_{m\in M}$ of ordinals, so that $\xi_m<\xi$ for $m\in M$ and
$\displaystyle{\sup_{m\in M}}\xi_m=\xi$. With suitable choices one can
define such systems that are useful for theoretical purposes or for
applications. In Theorem 1.6 we define a (Schreier type) system
$(\cca_{\xi})_{\xi<\omega_1}$ of uniform families, which in the
$\omega^a$--position for every $a<\omega_1$ has the family
$\ccb_a=\cca_{\omega^a}$ (Definition 1.5), a family similar to the
generalized Schreier set $\ccf_a$ (see Corollary 3.2) defined by Alspach
and Argyros in \cite{1}. Use of the system $(\ccf_a)_{a<\omega_1}$ has
proved fruitful, especially in connection with the theory of Banach spaces.
However, the system $(\ccf_a)_{a<\omega_1}$ is very difficult to employ in
inductive arguments owing mainly to lack of adequate interrelation of the
families $\ccf_a$, $a<\omega_1$, as there are missing families not defined
for all ordinals $\xi$ with $\omega^a<\xi<\omega^{a+1}$, $a<\omega_1$. The
introduction, in this paper, of the system $(\cca_{\xi})_{\xi<\omega_1}$
provides us with the correct amount of leeway to confront analogous
problems (see Section 3).

Our starting point is the following far--reaching generalization of the
classical Ramsey theorem.

\ \\
\indent
{\bf Theorem A} If $\ccf$ is a family of finite subsets of $\nn$, then for
every countable ordinal $\xi$, every infinite subset $M$ on $\nn$ and every
$\xi$--uniform family $\ccl$ on $M$ there exists an infinite subset $L$ of
$M$ such that either $\ccl\cap [L]^{<\omega}\subseteq \ccf$ or $\ccl\cap
[L]^{<\omega}\subseteq [M]^{<\omega}\smallsetminus \ccf$.

\ \\
\indent
A proof directly from the definitions involved is given in Theorem 2.2;
another proof, using the combinatorial theorems of Nash--Williams in
\cite{11} is given in \cite{13}.

For hereditary families of finite subsets of $\nn$ we prove a stronger
dichotomy result (Theorem 2.12, Th. B below). For the proof we introduce
the notion of the ``canonical representation'' for every finite subset of
$\nn$ with respect to a $\xi$--uniform family for every $\xi<\omega_1$
(Proposition 2.7).

\ \\
\indent
{\bf Theorem B} If $\ccf$ is a hereditary family of finite subsets of
$\nn$, then for every countable ordinal $\xi$, every infinite subset of
$\nn$ and every $\xi$--uniform family $\ccl$ on $M$ there exists an
infinite subset $L$ of $M$ such that either $\ccl\cap
[L]^{<\omega}\subseteq \ccf$, or $\ccf\cap [L]^{<\omega}\subseteq
\ccl^{*}\smallsetminus \ccl$ (where $\ccl^{*}$ is the family of
all the initial segments of the elements of $\ccl$).

\ \\
\indent
After that dichotomy result, with the help of the strong Cantor--Bendixson
index defined in \cite{2} and denoted by $s_M$, we describe when a
hereditary and pointwise closed family $\ccf$ of finite subsets of $\nn$ satisfies
one of the conditions given in Theorem B. A hereditary family $\ccf$ is
pointwise closed if and only if no infinite subset $M$ of $\nn$ exists such that
$[M]^{<\omega}\subseteq \ccf$ (Proposition 2.14). In fact, the following
result is proved in Theorem 2.16 and Remark 2.17.

\ \\
\indent
{\bf Theorem C} If $\ccf$ is a hereditary and pointwise closed family of finite
subsets of $\nn$, then for every countable ordinal $\xi$, every infinite
subset $M$ of $\nn$ and every $\xi$--uniform family $\ccl$ on $M$, the
following hold:

(i) If $s_M(\ccf)>\xi+1$ then there exists an infinite subset $L$ of $M$
such that
$$\ccl\cap [L]^{<\omega}\subseteq \ccf,$$

(ii) If $s_M(\ccf)<\xi+1$ then there exists an infinite subset $L$ of $M$
such that
$$\ccf\cap[L]^{<\omega}\subseteq \ccl^{*}\smallsetminus \ccl.$$

(iii) If $s_M(\ccf)=\xi+1$ then both alternatives ((i) and (ii)) may
materialize. 

\ \\
\indent
A consequence of these theorems is the existence, for a hereditary and
pointwise closed family $\ccf$ of finite subsets of $\nn$, of a
countable ordinal $\xi$ such that, for every system
$(\cca_{\zeta})_{\zeta<\omega_1}$ of uniform families the following obtain:
\begin{enumerate}
\item [(i)] For every $\zeta$ with $\zeta+1<\xi$ there exists an infinite
subset $L$ of $M$ such that
$$\cca_{\zeta}\cap[L]^{<\omega}\subseteq \ccf.$$ 
\item [(ii)] For every $\zeta$ with $\zeta<\zeta+1$ there exists an infinite
subset $L$ of $M$ such that
$$\ccf\cap[L]^{<\omega}\subseteq (\cca_{\zeta})^{*}\smallsetminus
\cca_{\zeta};$$ 
which gives that
$$\cca_{\zeta}\cap[L]^{<\omega}\subseteq [N]^{<\omega}\smallsetminus
\ccf.$$ 
\item [(iii)] If $\xi=\zeta+1$ both alternatives ((i) and (ii)) may materialize.
\end{enumerate}

Finally, for every hereditary family $\ccf$ of finite subsets of $\nn$
there exists an infinite subset $L$ of $\nn$ such that either
$[L]^{<\omega}\subseteq \ccf$, if $\ccf$ not closed, or
$[L]^{<\omega}\subseteq ([\nn]^{<\omega}\smallsetminus \ccf)_{*}$ (if
$\ccf$ is closed) (Corollary 2.15). $\ccf_{*}$ denotes the corresponding
hereditary family to $\ccf$).

\ \\
\indent
{\bf Notation and terminology:} We denote by $\nn$ the set of all natural
numbers. For an infinite subset $M$ of $\nn$ we denote by $[M]^{<\omega}$
the set of all finite subsets of $M$ and by $[M]$ the set of all infinite
subsets of $M$ (considering them as strictly increasing sequences).

If $s,t$ are finite subsets of $\nn$ then $s\preceq t$ means that $s$ is an
initial segment of $t$, while $s\prec t$ means that $s$ is a proper initial
segment of $t$. We write $s\leq t$ if $\max a\leq \min t$, while $s<t$ if
$\max s<\min t$.

Identifying every subset of $\nn$ with its characteristic function, we
topologize the set of all subsets of $\nn$ by the topology of pointwise
convergence. 

The generalized Schreier system $(\ccf_a)_{a<\omega_1}$, mentioned before,
has been defined in \cite{1} as follows:
$$\ccf_0=\{\{n\}:n\in \nn\};$$
if $\ccf_{\xi}$ has been defined
$$\ccf_{\xi+1}=\left\{\bigcup^n_{i=1}F_i:n\leq F_1<\ldots <F_n \ \mbox{ and
}\ F_i \in \ccf_{\xi}\right\};$$
if $\xi$ is limit choose $(\xi_n)_{n\in \nn}$ strictly increasing to $\xi$
and set
$$\ccf_{\xi}=\{F:F\in \ccf_{\xi_n} \ \mbox{ and }\ n\leq \min F\}.$$

\section{Systems of uniform families and Cantor--Bendixson index}

The definition of a uniform family (consisting of finite subsets of $\nn$),
stated below, 
is given by Pudl\'{a}k and R\"{o}dl in \cite{13}.
\begin{defi}
{\rm Let $M\in [\nn]$ and $\ccl$ be a family of finite subsets of $M$.
\begin{enumerate}
\item [(i)] For every $m\in M$, set $\ccl(m) =\{s\in
[M]^{<\omega}:\{m\}\cup s\in \ccl$ and $\{m\}<s\}$.
\item [(ii)] (Recursive definition of a uniform family)
\begin{enumerate}
\item [1.] $\ccl$ is {\bf 0--uniform on $M$} if $\ccl =\{\varnothing\}$;
\item [2.] if $\xi$ is a successor, countable ordinal, $\xi=\zeta+1$, then
$\ccl$ is $\xi${\bf--uniform on $M$} if $\varnothing \not\in \ccl$ and
the family $\ccl(m)$ is $\zeta$--uniform on $M\cap (m,+\infty)$ for every
$m\in M$; and 
\item [3.] if $\xi$ is a non--zero, limit countable ordinal then $\ccl$ is {\bf
$\xi$--uniform on $M$} if $\varnothing \not\in \ccl$ and there is an
increasing sequence $(\xi_m)_{m\in M}$ of ordinal numbers, smaller than
$\xi$, with $\displaystyle{\sup_m \xi_m=\xi}$ such that the family $\ccl
(m)$ is $\xi_m$--uniform on $M\cap (m,+\infty)$ for every $m\in M$.
\end{enumerate}
\item [(iii)] $\ccl$ is called {\bf uniform on $M$} if $\ccl$ is $\xi$--uniform
on $M$ for some countable ordinal $\xi$.
\item [(iv)] $\ccl^{*}=\{t\in [M]^{<\omega}:t$ is an initial segment of
some $s\in \ccl\}$.
\item [(v)] $\ccl_{*}=\{t\in [M]^{<\omega}:t\subseteq s$ for some $s\in
\ccl\}$. 
\item [(vi)] $\ccl$ is {\bf hereditary} if $\ccl_{*}=\ccl$.
\item [(vii)] $\ccl$ is {\bf Sperner} if there do not exist $s,t\in\ccl$
such that $s$ is a proper subset of $t$.
\item [(viii)] $\ccl$ is {\bf thin} if there do not exist $s,t\in\ccl$ such
that $s$ is a proper initial segment of $t$.
\end{enumerate}}
\end{defi}

Every family $\ccl \subset [M]^{<\omega}$ determines a partition
$\ccl=\displaystyle{\bigcup_{m\in M}\{\{m\}\cup s:m<s, s\in \ccl (m)\}}$;
and $\ccl$ is a $\xi$--uniform family precisely when  $\ccl(m)$ is
$\xi_m$--uniform on $M\cap (m,+\infty)$, with $\xi_m=\zeta$ for every
$m\in M$, if $\xi=\zeta+1$, and $(\xi_m)_{m\in M}$ an increasing to $\xi$
sequence if $\xi$ is limit. For example, if $\ccl$ is 1--uniform on $M$
then $\ccl=\{\{m\}:m\in M\}$ since $\ccl(m)=\{\varnothing\}$ (the only
0--uniform on $M\cap (m,+\infty)$).

Conversely, for every countable ordinal $\xi$ and $M\in [\nn]$ we can
construct a $\xi$--uniform family $\ccl$ on $M$ if we have for every $m\in
M$ a $\xi_m$--uniform family $\cca_m$ on $M\cap (m,+\infty)$,
($(\xi_m)_{m\in M}$ is as before). Indeed, the family
$\ccl=\displaystyle{\bigcup_{m\in M}\{\{m\} \cup s:s\in \cca_m\}}$ is
$\xi$--uniform on $M$, since $\ccl(m)=\cca_m$ for every $m\in M$.
\newtheorem{rem2}[defi]{Remarks}
\begin{rem2}
{\rm \begin{enumerate}
\item [(i)] For every $M\in [\nn]$ and every natural number $k$ there is
exactly one $k$--uniform family on $M$, namely the family $[M]^k$ of all
$k$--element subsets of $M$.
\item [(ii)] If $\ccl$ is a $\xi$--uniform family on $M$ $(M\in [\nn])$ and
$L\in [M]$, then, as can be proved by induction on $\xi$, $\ccl \cap
[L]^{<\omega}$ is $\xi$--uniform on $L$ (cf. \cite{13}).
\item [(iii)] Using (i) and (ii) we can describe a way of
constructing uniform families.

If $\ccl_{\xi}$ is a $\xi$--uniform family on $M$ (with $M\in [\nn]$) and
$k\in \nn$, then it is easy to see by induction on $k$ that the family
$$\ccl_{\xi+k}=\left\{ s\in [M]^{<\omega}:s=s_1\cup s_2 \mbox{ where
}s_1<s_2, s_1\in [M]^k \mbox{ and }s_2\in \ccl_{\xi}\right\}$$
is a $\xi+k$--uniform family on $M$.

If $\xi$ is a limit ordinal and $\ccl_{\beta}$ is a $\beta$--uniform family for
every $\beta<\xi$, then we choose an increasing sequence $(\xi_m)_{m\in
M}$ of ordinal numbers smaller than $\xi$ with
$\displaystyle{\sup_{m}\xi_m=\xi}$ and set 
$$\ccl_{\xi}=\left\{ s\in [M]^{<\omega}:s=\{m\}\cup s_1 \mbox{ where
}m\in M, \{m\}<s_1 \mbox{ and }s_1\in \ccl_{\xi_m}\right\};$$
the family $\ccl_{\xi}$ is in fact $\xi$--uniform on $M$, since
$(\ccl_{\xi})(m)=\ccl_{\xi_m}\cap [(m,+\infty)]^{<\omega}$ for every $m\in M$.
\item [(iv)] Every uniform family $\ccl$ on $M\in [\nn]$ is a maximal thin
subset of $[M]^{<\omega}$, (for the proof see \cite{13}). 
\item [(v)] A uniform family $\ccl$ on $M\in [\nn]$ is not necessarily
Sperner (see Example 1.12 below). 
However, 
for every uniform family $\ccl$ on $M$ there exists $L\in [M]$ such that
$\ccl\cap [L]^{<\omega}$ is Sperner (Corollary 2.4 below).
\end{enumerate}}
\end{rem2}

\ \\
\indent
Now, we will introduce the concept of a system of uniform families. A
system of uniform families on $M$ $(M\in [\nn])$ is a collection
$\cca=(\cca_{\xi})_{\xi<\omega_1}$, where each $\cca_{\xi}$ is
$\xi$--uniform on $M$, constructed in the way described in Remark 1.2
(iii), from uniform families $\cca_{\beta},\beta<\xi$, belonging to $\cca$.
The definition provides the necessary path, through which uniform families
are constructed and also gives the means of verification that a given
family is uniform.
\begin{defi}
{\rm Let $M\in [\nn]$ and $\cca_{\xi}\subseteq [M]^{<\omega}$ for every
countable ordinal $\xi$. The collection $\cca=(\cca_{\xi})_{\xi<\omega_1}$
is a {\bf system of uniform families on $M$} if:
\begin{enumerate}
\item [(i)] $\cca_{\xi}$ is $\xi$--uniform family on $M$ for every
$\xi<\omega_1$, and
\item [(ii)] For every $m\in M$ and $1\leq \xi<\omega_1$ 
$$\cca_{\xi}(m)=\cca_{\xi_m}\cap [(m,+\infty)]^{<\omega},$$
\end{enumerate}
where $\xi_m+1=\xi$ if $\xi$ is a successor ordinal; and $(\xi_m)_{m\in M}$
is an increasing sequence of ordinals smaller than $\xi$, with
$\displaystyle{\sup_m \xi_m=\xi}$, if $\xi$ is a limit ordinal.}
\end{defi}
\begin{rem2}
{\rm \begin{enumerate}
\item [(i)] If $(\cca_{\xi})_{\xi<\omega_1}$ is a system of uniform
families on $M$ and $L\in [M]$, then\\
$(\cca_{\xi}\cap [L]^{<\omega})_{\xi<\omega_1}$ is a system of uniform
families on $L$, according to Remark 1.2 (ii).
\item [(ii)] Every system of uniform families on $M$ is characterized by
the choices of the sequences $(\xi_m)_{m\in M}$ for every limit ordinal
$\xi$. Indeed, if for every limit
ordinal $\xi$ an increasing sequence $(\xi_m)_{m\in M}$ is given with
$\xi_m<\xi$ for every $m\in M$ and $\displaystyle{\sup_m\xi_m=\xi}$, then
we can define exactly one system of uniform families using these sequences
in the following way:
\begin{eqnarray*}
&& \cca_0=\{\varnothing\};\\
&&\cca_{\zeta+1}=\bigcup_{m\in M}\{\{m\}\cup s: s\in \cca_{\zeta}\cap
[(m,+\infty)]^{<\omega}\} \mbox{ for }\zeta <\omega_1; \mbox{ and}\\
&&\cca_{\xi}=\bigcup_{m\in M}\{\{m\}\cup s: s\in \cca_{\xi_m}\cap
[(m,+\infty)]^{<\omega}\},
\end{eqnarray*}
for $\xi$ limit, countable ordinal.
\end{enumerate}}
\end{rem2}

As we observed in Remark 1.2 (iii), for every $\xi$ with $\omega\leq
\xi<\omega_1$, there are continuum many $\xi$--uniform families. Indeed,
there are as many $\omega$--uniform families on $\nn$, as the multitude of
all the increasing, unbounded sequences of natural numbers.
Also, according to Remark 1.4 (ii) there are as many systems of uniform
families on $\nn$, as the multitude of all the choices of increasing
sequences $(\xi_n)_{n\in \nn}$, with $\xi_n<\xi$ for all $n\in \nn$ and
$\displaystyle{\sup_n \xi_n=\xi}$, for each countable limit ordinal $\xi$.

With suitable choices of sequences $(\xi_n)_{n\in \nn}$ one can define
interesting systems of uniform families. 
Below, in Theorem 1.6, we will define a Schreier type
system $\cca=(\cca_{\xi})_{\xi<\omega_1}$ of uniform families. This system
in the $\omega^{a}$--position has the uniform family
$\ccb_a=\cca_{\omega^a}$ (Definition 1.5 below) which is similar to the
Schreier set $\ccf_a$ (for every $a<\omega_1$) defined in \cite{1}.
\begin{defi}
{\rm (Schreier type system of uniform families) \\
(1) We define inductively for every
$a<\omega_1$ the families $\ccb_a\subseteq [\nn]^{<\omega}$ as follows:
\begin{enumerate}
\item [(i)] $\ccb_0=\{\{n\}:n\in \nn\}$;
\item [(ii)] If the family $\ccb_a$ has been defined, let
$$\ccb_{a+1}=\{s\subseteq \nn: s=\bigcup^n_{i=1}s_i \ \mbox{ where }n=\min
s_1,s_1<\ldots<s_n \mbox{ and }s_1,\ldots,s_n\in \ccb_a\};\mbox{ and}$$
\item [(iii)] If $a$ is a limit countable ordinal and the families
$\ccb_{\zeta}$ have been defined for each $\zeta<a$, let
$$\ccb_a=\{s\subseteq \nn:s\in \ccb_{a_n} \mbox{ with }n=\min s\},$$
where $(a_n)$ is a fixed increasing sequence of ordinal numbers smaller
than $a$ with $\displaystyle{\sup_n a_n=a}$.
\end{enumerate}}
\end{defi}

\noindent
(2) We set $\cca_{\omega^a}=\ccb_a$ for all ordinals $a<\omega_1$, and we
complete the system of uniform families as follows:
\begin{enumerate}
\item [(i)] $\cca_0=\{\varnothing\}$
\item [(ii)] if $\xi<\omega_1$, and the family $\cca_{\xi}$ has been
defined, then set 
$$\cca_{\xi+1}=\{s\subseteq \nn: s=\{n\}\cup s_1 \mbox{ where }n\in \nn,
\{n\} <s_1 \mbox{ and }s_1\in \cca_{\xi}\};\mbox{ and}$$
\item [(iii)] if $\xi$ is a limit countable ordinal and the families
$\cca_{\zeta}$ have been defined for every $\zeta<\xi$ and if $\xi$ has the
form $\xi=\displaystyle{\sum^m_{i=1}p_i\omega^{a_i}}$, where
$m,p_1,\ldots,p_m\in \nn$ and $a_1>\ldots>a_m>0$ are ordinal numbers, then
we set 
\begin{eqnarray*}
\cca_{\xi} &=& \{s\subseteq \nn: s=\bigcup^m_{i=1}s_i\mbox{ where
}s_m<\ldots<s_1,s_i=F_1^i\cup \ldots\cup F^i_{p_i} \mbox{ with
}F_1^i<\ldots<F^i_{p_i}\\ 
&& \mbox{ and}F_1^i,\ldots,F^i_{p_i}\in \ccb_{a_i}\mbox{ for every }1\leq i \leq m\}.
\end{eqnarray*}
\end{enumerate}
\newtheorem{theo}[defi]{Theorem}
\begin{theo}
The collection $(\cca_{\xi})_{\xi<\omega_1}$ is a system of uniform
families on $\nn$.
\end{theo}
\pr $\cca_0=\{\varnothing\}$, so it is $0$--uniform on $\nn$. We assume
that for every $\zeta<\xi$ 
the families $\cca_{\zeta}$ are $\zeta$--uniform on
$\nn$ and also that $\cca_{\zeta}(n)=\cca_{\zeta_n}\cap
[(n,+\infty)]^{<\omega}$ for every $n\in \nn$, where $\zeta_n+1=\zeta$ for
every $n\in \nn$, if $\zeta$ is a successor ordinal; and $(\zeta_n)$ is an
increasing sequence of ordinals smaller than $\zeta$ with $\sup
\zeta_n=\zeta$, if $\zeta$ is limit.

Let $\xi=\zeta+1$ be a successor ordinal. According to the definition of
$\cca_{\zeta+1}$ (1.5) 
$\cca_{\xi}(n)=\cca_{\zeta}\cap [(n,+\infty)]^{<\omega}$ for every $n\in
\nn$. Hence, $\cca_{\xi}$ is $\xi$--uniform on $\nn$, since $\cca_{\zeta}$ is
$\zeta$--uniform (Remark 1.2 (ii)).

Let $\xi$ be a limit ordinal. We will check all particular cases:

\noindent
(1) If $\xi=\omega$ then for every $n\in \nn$
$$\cca_{\omega}(n)=\ccb_1(n)=\{s:\{n\}\cup s\in \ccb_1 \mbox{ and }\{n\}
<s\}= [(n,+\infty)]^{n-1}=\cca_{n-1}\cap [(n,+\infty)]^{<\omega}.$$
Hence $\cca_{\omega}$ is $\omega$--uniform on $\nn$ and $\omega_n=n-1$ for every
$n\in \nn$. 

\noindent
(2) If $\xi=\omega^{a+1}$ then for every $n\in \nn$
\begin{eqnarray*}
\cca_{\omega^{a+1}}(n) &=& \ccb_{a+1}(n)=\{s: s=s_1\cup s_2 \mbox{ where
}\{n\} <s_1<s_2,s_1\in \ccb_a(n) \mbox{ and}s_2\in \cca_{(n-1)\omega^a}\}=\\
&& \cca_{(n-1)\omega^a+(\omega^a)_n}\cap [(n,+\infty)]^{<\omega}.
\end{eqnarray*}
Hence $\cca_{\omega^{a+1}}$ is $\omega^{a+1}$--uniform on $\nn$ and
$(\omega^{a+1})_n=(n-1)\omega^a+(\omega^a)_n$ for every $n\in \nn$.

\noindent
(3) If $\xi=\omega^a$ for a limit ordinal $a$ then
$\cca_{\omega^a}=\ccb_a$. 
Let $(a_n)$ be the fixed sequence of ordinal numbers which is used in the
definition of $\ccb_a$ (Definition 1.5). For every $n\in \nn$ we have 
$$\cca_{\omega^a}(n)=\cca_{(\omega^{a_n})_n}\cap
[(n,+\infty)]^{<\omega}.$$ 
Hence, $\cca_{\omega^a}$ is $\omega^a$--uniform on $\nn$ and
$(\omega^a)_n=(\omega^{a_n})_n$ for every $n\in \nn$.

\noindent
(4) If $\xi =p\omega^a$, where $p\in \nn$ and $0<a<\omega_1$, then for every
$n\in \nn$
\begin{eqnarray*}
\cca_{p\omega^{a}}(n) &=& \{s: s=s_1\cup s_2 \mbox{ where
}\{n\} <s_1<s_2,s_1\in \cca_{(\omega^a)_n}\mbox{ and}s_2\in \cca_{(p-1)\omega^a}\}=\\
&& \cca_{(p-1)\omega^a+(\omega^a)_n}\cap [(n,+\infty)]^{<\omega}.
\end{eqnarray*}                                                             
Hence, $\cca_{p\omega^a}$ is $p\omega^a$--uniform on $\nn$ and
$(p\omega^a)_n=(p-1)\omega^a+(\omega^a)_n$ for every $n\in \nn$.

\noindent
(5) Finally, if $\xi=\displaystyle{\sum^m_{i=1}p_i\omega^{a_i}}$, where
$m,p_1,\ldots,p_m\in \nn$ and $a_1>\ldots>a_m>0$, then for every $n\in \nn$
\begin{eqnarray*}
\cca_{\xi}(n) &=& \{s: s=s_1\cup s_2 \mbox{ where
}\{n\} <s_1<s_2,s_1\in \cca_{(p_m \omega^{a_m})_n}\mbox{ and}s_2\in 
\cca_{\beta}\}=\\
&=&\cca_{\beta+(p_m\omega^{a_m})_n}\cap [(n,+\infty)]^{<\omega},
\end{eqnarray*}
where $\beta=\displaystyle{\sum^{m-1}_{i=1}p_i\omega^{a_i}}$.\\
Hence, $\cca_{\xi}$ is $\xi$--uniform on $\nn$ and
$\xi_n=\displaystyle{\sum^{m-1}_{i=1}p_i\omega^{a_i}+(p_m\omega^{a_m})_n}$
for every $n\in \nn$.

This completes the proof of the theorem.
\newtheorem{coro}[defi]{Corollary}
\begin{coro}
{\rm For every $M\in [\nn]$, the collection
$(\cca^M_{\xi})_{\xi<\omega_1}$, where
$\cca^M_{\xi}=\cca_{\xi}\cap[M]^{<\omega}$ for every $\xi<\omega_1$, is a
system of uniform families on $M$.}
\end{coro}
\pr This is a consequence of Theorem 1.6 and Remark 1.4 (i).

\ \\
\indent
It would be very complicated to prove directly that the family $\ccb_a$ is
$\omega^a$-- uniform on $\nn$ for every $a<\omega_1$, but using the notion
of a system of uniform families the proof is immediate after Theorem 1.6.
\begin{coro}
{\rm For every countable ordinal $a$ and $M\in [\nn]$ the family
$\ccb_a\cap [M]^{<\omega}$ is $\omega^a$--uniform on $M$.}
\end{coro}
\pr We have $\ccb_a=\cca_{\omega^a}$ for every $a<\omega_1$.
\begin{rem2}
{\rm \begin{enumerate}
\item [(i)] In the definition of the Schreier type system
$\cca=(\cca_{\xi})_{\xi<\omega_1}$ of uniform families we have choices of
increasing sequences $(\xi_n)_{n\in \nn}$ (in fact $\xi_n=\omega^{a_n}$ for
every $n\in \nn$) with $\xi_n<\xi$ for all $n\in \nn$ and
$\displaystyle{\sup_n\xi_n=\xi}$ only in the cases $\xi=\omega^a$, where $a$
is a limit countable ordinal. In the other limit countable ordinals $\xi$ we
use concrete sequences depending on $\xi$ and the previous choices.
\item [(ii)] It is easy to see that $\ccb_a\subseteq \ccf_a$ for every
$a<\omega_1$. In general, the hereditary family $(\ccb_a)_{*}$ of all the
subsets of the elements of $\ccb_a$ is not equal to $\ccf_a$. However, in
Section 3 (Proposition 3.1) we will prove that for every $M\in [\nn]$
there exists $L=(\ell_n)_{n\in \nn}\in [M]$ such that $\ccf_a(L)\subseteq
(\ccb_a)_{*}\cap [M]^{<\omega}$, where
$$\ccf_a(L)=\{(\ell_{n_1},\ldots,\ell_{n_k})\subseteq L:(n_1,\ldots,n_k)\in
\ccf_a\}.$$
\end{enumerate}}
\end{rem2}

At this point the reader might think that the definition of a system of
uniform families is unneccesarily cumbersome. It bears similarity to the
various Schreier--type system $(\ccf_a)_{a<\omega_1}$ used in the
literature (e.g. Alspach--Argyros (\cite{1}),
Argyros--Mercourakis--Tsarpalias (\cite{2}), Farmaki (\cite{6}, \cite{7}),
Kyriakouli--Negrepontis (\cite{10}), Odell--Tomczak--Wagner (\cite{12}) and
others). However, the system 
$(\ccf_a)_{a<\omega_1}$ is very difficult to employ in inductive arguments,
owing on the one hand to the concrete and fixed nature of the definition of
$\ccf_a,a<\omega_1$ and also, and more significantly on a rather more
hidden aspect of their interrelation (for different ordinals). We can
clarify the precise relation between the system $(\ccf_a)_{a<\omega_1}$ and
the system of uniform families $(\cca_{\xi})_{\xi<\omega_1}$, if we think
that each family $\ccf_a$ is related not to the family $\cca_a$ but to the
uniform family $\cca_{\omega^a}$. In other words, the difficulty in
employing Schreier--type systems in inductive arguments lies with the fact
that, e.g. there are missing families, not defined for all ordinals $\xi$
with $\omega^a<\xi<\omega^{a+1},a<\omega_1$. This filling up of the
intermediate gaps was in effect performed in a special case, arising in
Banach space theory, in our earlier work in \cite{7}.

Thus, returning to the difficulty in employing induction on the Schreier
sets $\cca_a,a<\omega_1$ owing to their fixed nature, it will be seen
clearly in Section 2 below that the notion of a uniform family ($\cca_{\xi}$
is uniform for every $\xi<\omega_1$) provides us with the correct amount of
leeway, $a$ leeway that is precisely missing from the system
$(\ccf_a)_{a<\omega_1}$. 

\ \\
\indent
In the following we will estimate the strong Cantor--Bendixson index of a
uniform family.
This index (see Definition 1.10 below) is analogous to the well--known
Cantor--Bendixson index (\cite{4}, \cite{5}) and has been defined in \cite{2}. Here, we will use
a different notation in order to avoid some misinterpretations.

We will prove in Proposition 1.18 below that, for every $\xi<\omega_1$, $M\in
[\nn]$,  the corresponding
hereditary family $\ccl_{*}$ of a $\xi$--uniform family $\ccl$ on $M$ has
strong Cantor--Bendixson index 
on $M$ equal to $\xi+1$. Hence, if $(\cca_{\xi})_{\xi<\omega_1}$ is a system
of uniform families, then the collection $((\cca_{\xi})_{*})_{\xi<\omega_1}$
contains hereditary families of arbitrary index.
\begin{defi}
{\rm (\cite{2}) Let $\ccf$ be a hereditary and pointwise closed family of
finite subsets of $\nn$. For $M\in [\nn]$ we define the {\bf strong
Cantor
-- Bendixson derivatives} $(\ccf)^{\xi}_M$ of $\ccf$ on $M$ for every
$\xi<\omega_1$ as follows:
$$(\ccf)^1_M=\{A\in \ccf [M]: A \mbox{ is a cluster point of }\ccf [A\cup L]
\mbox{ for each }L\in [M]\};$$
(where, $\ccf[M]=\ccf\cap [M]^{<\omega}$).\\
If $(\ccf)^{\xi}_M$ has been defined, then
$$(\ccf)_M^{\xi+1}=((\ccf)^{\xi}_M)^1_M.$$
If $\xi$ is a limit ordinal and $(\ccf)^{\beta}_M$ have been defined for
each $\beta<\xi$, then 
$$(\ccf)^{\xi}_M=\displaystyle{\bigcap_{\beta<\xi}}(\ccf)^{\beta}_M.$$}
\end{defi}

The {\bf strong Cantor -- Bendixson index of $\ccf$ on $M$} is defined to
be the smallest countable ordinal $\xi$ such that $(\ccf)^{\xi}_M=\varnothing$. We
denote this index by $s_M(\ccf)$.

We define the {\bf strong Cantor -- Bendixson index} $s(\ccf)$ of $\ccf$ to
be $s(\ccf)=s_M(\ccf)$, where $M=\{n\in \nn:\{n\}\in \ccf\}$ is the support
of $\ccf$.
\begin{rem2}
{\rm \begin{enumerate}
\item [(i)] Of course, the strong Cantor -- Bendixson index is a successor
ordinal.
\item [(ii)] If $\ccf_1,\ccf_2 \subseteq [\nn]^{<\omega}$ are hereditary
and closed and $\ccf_1\subseteq \ccf_2$ then $s_M(\ccf_1)\leq s_M (\ccf_2)$
for every $M\in [\nn]$.
\item [(iii)] $s_M(\ccf)=s_M(\ccf[M])$ for every $M\in [\nn]$.
\item [(iv)] For every $M\in [\nn]$ and $A \in [M]^{<\omega}$ according to
a remark in \cite{9} we have:\\
$A\in (\ccf)^1_M$ if and only if the set $\{m\in M: A\cup
\{m\}\not\in\ccf\}$ is finite. 
\item [(v)] Using the previous remark (iv) can be proved by induction that
for every $L\in [M]$ and $\beta<\omega_1$
$$A\cap L\in (\ccf)^{\beta}_L \mbox{ if }A\in (\ccf)^{\beta}_M.$$
Hence, $s_L(\ccf)\geq s_M(\ccf)$. (see also \cite{2}).
\item [(vi)] If $L$ is almost contained in $M$, then 
$$s_L(\ccf)\geq s_M(\ccf).$$
\item [(vii)] For every $a<\omega_1$, let $\ccf_a$ be the Shreier family.
Then for every $M\in [\nn]$
$$s_M(\ccf_a)=\omega^a+1.$$
(see \cite{2}).
\end{enumerate}}
\end{rem2}

In the following we will give the precise relation between the strong
Cantor--Bendixson derivatives of the corresponding hereditary family $\ccl_{*}$
of a given family $\ccl\subseteq [\nn]^{<\omega}$ and the derivatives of
the families $(\ccl(n))_{*}$ for every $n\in \nn$. After that, we will
calculate the strong Cantor--Bendixson index of a uniform family.

First of all we must notice that the families $(\ccl (n))_{*}$ and
$\ccl_{*}(n)$ are in general different as we can see from the following
example. 
\newtheorem{exam}[defi]{Example}
\begin{exam}
{\rm For every $n\in \nn$ choose the following member of the Schreier type
system $(\cca_{\xi})_{\xi<\omega_1}$ (Definition 1.5).

$\ccl_1=[\nn]^5=\cca_5,$

$\ccl_2=\cca_{\omega}$, and

$\ccl_n=\cca_{\omega+n}$ for every $n>2$.

Set
$$\ccl=\bigcup_{n\in \nn}\{\{n\} \cup s: s\in \ccl_n\mbox{ and }\{n\}<s\}.$$ 
Then $\ccl$ is a $2\omega$--uniform family. Let $s=(2,3,4,5,6)$ and
$t=(1,2,3,4,5,6)$. Since $t\in \ccl$ we have that $s\in \ccl_{*}$ and
consequently that $s_1=(3,4,5,6)\in 
\ccl_{*}(2)$. As we can see, $\ccl (2) =\cca_{\omega}\cap
[(2,+\infty)]^{<\omega}$ and of course $s_1\not\in (\ccl(2))_{*}$.

It is remarkable that the family $\ccl$ is not Sperner (see Remark 1.2
(v)), since $F=(2,3,4,5)\in \ccl$ and $F\subseteq t$, $(t\in\ccl)$.}
\end{exam}
\newtheorem{lemm}[defi]{Lemma}
\begin{lemm}
{\rm Let $\beta$ be a countable ordinal and $\ccl\subseteq [\nn]^{<\omega}$
such that $\ccl_{*}$ and $\ccl(n)_{*}$ are closed for every $n\in \nn$. If
$A\in (\ccl (n)_{*})^{\beta}_M$ for some $n\in \nn$ and $M\in [\nn]$, then
$\{n\} \cup A\in (\ccl_{*})^{\beta}_M$.}
\end{lemm}
\pr We use induction on $\beta$. Let $A\in (\ccl(n)_{*})^1_M$
for some $n\in \nn$ and $M\in [\nn]$. Since
$$\{m\in M:A\cup \{m\}\in \ccl(n)_{*}\}\subseteq \{m\in M:A\cup \{m\}\cup
\{n\} \in \ccl_{*}\},$$
we have, according to Remark 1.11 (iv), that $A\cup \{n\}\in
(\ccl_{*})^{\beta}_M$. 

Suppose that the assertion holds for all ordinals $\zeta$ with
$\zeta<\beta$. If $A\in (\ccl(n)_{*})^{\zeta+1}_M$ then $\{n\}\cup A\in
(\ccl_{*})^{\zeta+1}_M$, since
$$\{m\in M:A\cup \{m\}\in (\ccl(n)_{*})_M^{\zeta}\}\subseteq \{m\in M:A\cup 
\{m\}\cup\{n\} \in (\ccl_{*})^{\zeta}_M\};$$
according to the induction hypothesis.

The case where $\beta$ is limit ordinal is trivial.
\newtheorem{prop}[defi]{Proposition}
\begin{prop}
{\rm Let $\ccl\subseteq [\nn]^{<\omega}$ such that $\ccl_{*}$ and
$\ccl(n)_{*}$  are closed for every $n\in \nn$ and $M\in [\nn]$.
\begin{enumerate}
\item [(i)] If there exists $L\in [M]$ such that
$s_M(\ccl(n)_{*})=\xi$ for every $n\in L$ then $s_L(\ccl_{*})\geq \xi+1$.
\item [(ii)] Let $\xi_n=s_M(\ccl(n)_{*})$ for every $n\in \nn$. If there
exists $L\in [M]$ such that $\xi_n<\xi=\displaystyle{\sup_{n\in L}}\xi_n$
for every $n\in L$, then $s_M(\ccl_{*})\geq \xi+1$.
\end{enumerate}}
\end{prop}
\pr \begin{enumerate}
\item [(i)] Let $\xi=\beta+1$ and $L\in [M]$ such that
$s_M(\ccl(n)_{*})=\xi$ for every $n\in L$. Then $\varnothing \in
(\ccl(n)_{*})^{\beta}_M$ for every $n\in L$. According to Lemma 1.13 we
have $\{n\}\in (\ccl_{*})^{\beta}_M$ for every $n\in L$ and then $\{n\}\in
(\ccl_{*})^{\beta}_L$ for every $n\in L$ by Remark 
1.11 (v). From Remark
1.11 (iv) we have $\varnothing \in (\ccl_{*})^{\beta+1}_L$ and therefore
$s_L(\ccl_{*})\geq \xi+1$.
\item [(ii)] In this case $\xi$ is a limit ordinal. Since the $\xi_n$ are
successor ordinals we set $\xi_n=\beta_n+1$ for every $n\in \nn$. According
to our hypothesis $\varnothing\in (\ccl(n)_{*})^{\beta_n}_M$ for every $n\in L$.
Hence $\varnothing\in (\ccl_{*})^{\beta_n}_M$ for every $n\in L$. Since
$\displaystyle{\sup_{n\in L}\beta_n=\xi}$ and $\beta_n<\xi$ we have
$\varnothing\in (\ccl_{*})^{\xi}_M$ and therefore $s_M(\ccl_{*})\geq \xi+1$.
\end{enumerate}
\begin{lemm}
{\rm Let $\beta$ be a countable ordinal and $\ccl\subseteq [\nn]^{<\omega}$
such that $\ccl_{*}$ and $\ccl(n)_{*}$ are closed for every $n\in \nn$. If
$A\neq \varnothing$ and $A\in (\ccl_{*})^{\beta}_M$ for some $M\in [\nn]$, then
there exist $\ell\in\nn$ with $\ell\leq \min A$ and $L\in [M]$ such that 
$A\smallsetminus \{\ell\}\in (\ccl(\ell)_{*})^{\beta}_L$.}
\end{lemm}
\pr We use induction on $\beta$. Let $A\neq \varnothing$ and $A\in
(\ccl_{*})^1_M$. According to Remark 1.11 (iv) the set $M_A=\{m\in M:A\cup
\{m\}\in \ccl_{*}$ and $\min A\leq m\}$ is almost equal to $M$. For each
$m\in M_A$ there exists $s_m\in \ccl$ such that $A\cup \{m\}\subseteq s_m$.
Set
$$\ell=\min \{n\in \nn:\mbox{ the set }\{m\in M_A:\min s_m=n\} \mbox{ is
infinite }\}.$$
Of course $\ell\leq A$. Set $L=\{m\in M_A:\min s_m=\ell\}\cup A$. Then
$L\in [M]$ and $A\smallsetminus \{\ell\}\in (\ccl (\ell)_{*})^1_L$, as
required. 

Suppose now that the assertion holds for all ordinals $\zeta$ with
$\zeta<\beta$ and let $\beta=\zeta+1$. If $A\neq \varnothing$ and $A\in
(\ccl_{*})^{\zeta+1}_M$ then according to Remark 1.11 (iv) the set\\
$M_A=\{m\in M:A\cup \{m\}\in (\ccl_{*})^{\zeta}_M$ and $\min A\leq m\}$ is
almost equal to $M$. Let $m_1=\min M_A$. By the induction
hypothesis there exist $\ell_1\in \nn$ and $L_1\in [M_A]$ with $\ell_1\leq
\min A$ such that $A\cup \{m_1\}\smallsetminus \{\ell_1\}\in (\ccl
(\ell_1)_{*})^{\zeta}_{L_1\cup A}$, since $A\cup \{m_1\}\in
(\ccl_{*})^{\zeta}_{M_A\cup A}$ (Remark 1.11 (v)). Let $m_2\in L_1$ and
$m_2>m_1$. Since $A\cup \{m_2\}\in (\ccl_{*})^{\zeta}_{L_1\cup A}$ there
exists $\ell_2\in \nn$ with $\ell_2\leq \min A$ and $L_2\in [L_1]$ such
that $A\cup \{m_2\}\smallsetminus \{\ell_2\}\in (\ccl
(\ell_2)_{*})^{\zeta}_{L_2\cup A}$. We continue analogously setting $m_3\in
L_2$ with $m_3>m_2$ and so on.

Hence we construct an increasing sequence $(m_i)^{\infty}_{i=1}$ in $M_A$, a
sequence $(\ell_i)^{\infty}_{i=1}$ in $\nn$, with \\
$1\leq \ell_i\leq \min A$ for every $i\in \nn$, and a decreasing sequence
$(L_i)_{i=1}^{\infty}$ in $[M_A]$ such that 
$$A\cup \{m_i\}\smallsetminus \{\ell_i\}\in (\ccl
(\ell_i)_{*})^{\zeta}_{L_i\cup A}$$
for every $i\in \nn$.

We can find $\ell \in \nn$ with $1\leq \ell \leq \min A$ such that the set
$I=\{i\in \nn: \ell_i=\ell\}$ is infinite. Set $L=\{m_i:i\in I\}\cup A$.
Then $A\smallsetminus \{\ell\}\in (\ccl(\ell)_{*})^{\zeta+1}_L$, as
required. 

In the case where $\beta$ is a limit ordinal and $A\in
(\ccl_{*})_M^{\beta}$, $A\neq \varnothing$, we fix a strictly increasing sequence
$(\zeta_i)^{\infty}_{i=1}$ of ordinals with $\zeta_i<\beta$ for every $i\in
\nn$ and $\displaystyle{\sup_i \zeta_i=\beta}$. Then $A\in
(\ccl_{*})^{\zeta_i}_M$ for every $i\in \nn$.
According to the induction hypothesis if $M_A=\{m\in M:\min A\leq m\}$
there exist $\ell_1\in\nn$ with $\ell_1\leq \min A$ and $L_1\in [M_A]$ such
that $A\smallsetminus \{ \ell_1\}\in (\ccl (\ell_1)_{*})^{\zeta_1}_{L_1\cup
A}$.

Since $A\in (\ccl_{*})^{\zeta_2}_{L_1\cup A}$ there exists $\ell_2\in \nn$
with $\ell_2\leq \min A$ and $L_2\in [L_1]$ such that $L_2\neq L_1$ and 
$$A\smallsetminus \{ \ell_2\}\in (\ccl (\ell_2)_{*})^{\zeta_2}_{L_2\cup A}$$
So we construct a sequence $(\ell_i)^{\infty}_{i=1}$ with $1\leq \ell_i\leq
\min A$ and a strictly decreasing sequence $(L_i)^{\infty}_{i=1}$ in
$[M_A]$ such that
$$A\smallsetminus \{ \ell_i\}\in (\ccl (\ell_i)_{*})^{\zeta_i}_{L_i\cup A}$$
for every $i\in \nn$.

We can find $\ell$ with $1\leq \ell\leq \min A$ such that the set $I=\{i\in
\nn:\ell_1=\ell\}$ is infinite. Set $L=\{\min L_i:i\in I\}\cup A$.
Then $A\smallsetminus \{ \ell\}\in (\ccl (\ell)_{*})^{\zeta_i}_{L}$ for
every $i\in I$. Since $\displaystyle{\sup_{i\in I}\zeta_i=\beta}$, we have
that 
$$A\smallsetminus \{ \ell\}\in (\ccl (\ell)_{*})^{\beta}_{L}.$$ 
This completes the proof.
\begin{prop}
{\rm Let $\ccl \subseteq [\nn]^{<\omega}$ such that $\ccl_{*}$ and $\ccl
(n)_{*}$ are closed for every $n\in \nn$ and $M\in [\nn]$. If $\xi=\sup
\{s_L(\ccl (n)_{*}):n\in \nn$ and $L\in [M]\}$ then $(\ccl_{*})^{\xi}_M
\subseteq \{\varnothing\}$ and therefore $s_M(\ccl_{*})\leq \xi+1$.}
\end{prop}
\pr Let $A\in (\ccl_{*})^{\xi}_M$ and $A\neq \varnothing$. According to lemma 1.15
there exist $n\in \nn$ and $L\in [M]$ such that $A\smallsetminus \{n\}\in
(\ccl (n)_{*})^{\xi}_L$. Hence, $s_L(\ccl(n)_{*}) \geq \xi+1$. A
contradiction, which finishes the proof.

\ \\
\indent
After Propositions 1.14 and 1.16 we will see in Theorem 1.18 that the
definition of a uniform family is the most suitable and least complicated in
order to ensure that every $\xi$--uniform family on $M$ $(M\in [\nn])$ is
thin (this arises from the condition $\varnothing\not\in \ccl$ for every
$\zeta$--uniform family with $1\leq \zeta$ in Definition 1.1) and the
corresponding hereditary family has strong Cantor--Bendixson index on $L$
equal to $\xi+1$, for every $L\in [M]$.
\begin{lemm}
{\rm Let $M\in [\nn]$ and $\ccl$ a $\xi$--uniform family on $M$, for some
$\xi<\omega_1$. Then $\ccl_{*}$ is closed.}
\end{lemm}
\pr This is easily proved by induction on $\xi$.
\begin{theo}
Let $\xi$ be a countable ordinal, $M\in [\nn]$ and $\ccl$ a $\xi$--uniform
family 
on $M$. Then for every $L\in [M]$ we have
$(\ccl_{*})^{\xi}_L=\{\varnothing\}$ and $s_L(\ccl_{*})=\xi+1$.
\end{theo}
\pr We use induction on $\xi$. Let $\xi=1$. Then
$\ccl_{*}=\{\{n\}:n\in M\}\cup \{\varnothing\}$. Hence
$(\ccl_{*})^1_L=\{\varnothing\}$ 
and therefore $s_L(\ccl_{*})=2$ for every $L\in [M]$.

Suppose the assertion holds for every ordinal number $\beta$ with $\beta<\xi$.
In case $\ccl$ is a $\zeta+1$--uniform on $M$ the families $\ccl(n)$
are $\zeta$--uniform on $M_n=M\cap (n,+\infty)$ for every $n\in M$. Hence
according to the induction hypothesis $s_L(\ccl(n)_{*})=\zeta+1=\xi$ for
every $L\in [M]$ (cf. Remark 1.11 (vi)). By Proposition 1.14 (i)
we have $s_L(\ccl_{*})\geq \xi+1$ for every $L\in [M]$. Hence
$(\ccl_{*})^{\xi}_L\neq \varnothing$ for every $L\in [M]$. On the other
hand, according to Proposition 1.16 we have $(\ccl_{*})^{\xi}_L\subseteq
\{\varnothing\}$ for every $L\in [M]$. Hence $(\ccl_{*})^{\xi}_L=\{\varnothing\}$ and
therefore $s_L(\ccl_{*})=\xi+1$ for every $L\in [M]$.

In the case where $\ccl$ is a $\xi$--uniform family on $M$ for a limit ordinal
$\xi$ we have that $\ccl(n)$ are $\beta_n$--uniform on $M\cap (n,+\infty)$
for every $n\in M$, where $(\beta_n)$ is a sequence of ordinals smaller
than $\xi$ with $\displaystyle{\sup_{n\in M}\beta_n=\xi}$.
According to the induction hypothesis and Remark 1.11 (vi) we have
$s_L(\ccl(n)_{*})=\beta_n+1=\xi_n$ for every $L\in [M]$ and $n\in M$. From
Proposition 1.14 (ii) we have that $s_L(\ccl_{*})\geq \xi+1$ for every
$L\in [M]$ hence $(\ccl_{*})^{\xi}_L\neq \varnothing$ for every $L\in [M]$. On the
other hand $(\ccl_{*})^{\xi}_L\subseteq \{\varnothing\}$ for every $L\in [M]$
according to Proposition 1.16. Hence $(\ccl_{*})^{\xi}_L=\{\varnothing\}$ and
therefore $s_L(\ccl_{*})=\xi+1$ for every $L\in [M]$.

The proof is complete.
\begin{coro}
{\rm Let $M\in [\nn]$ and $(\cca^M_{\xi})_{\xi<\omega_1}$ a system of
uniform families on $M$. Then \\
$s_L\left( (\cca_{\xi}^M)_{*}\right)=\xi+1$
for every $\xi<\omega_1$ and $L\in [M]$.}
\end{coro}

Hence a system of uniform families on $M$ is an appropriate selection
$(\cca_{\xi})_{\xi<\omega_1}$ of
thin subfamilies of $[M]^{<\omega}$ with strong Cantor--Bendixson index on
$L$ each countable ordinal, for every $L\in [M]$.

\section{Ramsey dichotomies with ordinal index}

We start this section with an equivalent formulation of the classical
Ramsey theorem (\cite{14}).
\begin{theo}
(Ramsey) For any positive integers $r$ and $k$ if we partition the family
$[M]^{<\omega}$ of all the finite subsets of an infinite set $M$ into
$k$--parts, then there is an infinite subset $L$ of $M$, all $r$--tuples of
which belong to the same class of the partition.
\end{theo}

In the following we will show how the concept of $\xi$--uniform families
can be applied to provide a far--reaching generalization of the classical
Ramsey theorem. This happens because general $\xi$--uniform families share
with the family $[M]^r$ of all $r$--tuples of $M$ occuring in the Ramsey
theorem the following properties: (a) they are thin and (b) the
Cantor--Bendixson index does not dicrease when we restrict ourselves to any
infinite subset of $M$.

Our proof will be an elementary one, directly from the definitions
involved. Another proof can be obtained, using the combinatorial theorem
of Nash--Williams in \cite{11} (see also \cite{13}).
\begin{theo}
Let $M$ be an infinite subset of $\nn$, $\{P_1,P_2\}$ a partition of the
set $[M]^{<\omega}$ of all finite subsets of $M$, $\xi$ a countable ordinal
number and $\ccl$ a $\xi$--uniform family on $M$. Then there exists an
infinite subset $L$ of $M$ such that either $\ccl\cap
[L]^{<\omega}\subseteq P_1$ or $\ccl\cap [L]^{<\omega}\subseteq P_2$.
\end{theo}
\pr We will prove the theorem by induction on $\xi$.

Let $\xi=1$. Then $\ccl=\{\{m\}:m\in M\}$. Set
\begin{eqnarray*}
M_1&=&\{m\in M:\{m\}\in P_1\} \mbox{ and}\\
M_2&=&\{m\in M:\{m\}\in P_2\}.
\end{eqnarray*}
If $M_1$ is infinite then the theorem holds for $L=M_1$, otherwise it holds
for $L=M_2$.

Assume that the theorem is valid for every ordinal $\zeta$ with $\zeta<\xi$
and let $\ccl$ be a $\xi$--uniform family on $M$. Then, according to
Definition 1.1, there exists a sequence $(\xi_m)_{m\in M}$ of ordinal
numbers such that $\xi_m<\xi$ for every $m\in M$ and the family $\ccl(m)$
is $\xi_m$--uniform on\\
$M\cap (m,+\infty)$.

Let $m_1=\min M$ and $M_1=M\cap (m_1,+\infty)$. Set 
\begin{eqnarray*}
P^1_1 &=& \{s\subseteq M:\{m_1\}\cup s\in P_1 \mbox{ and }\{m_1\}<s\}
\mbox{ and}\\
P^1_2 &=& \{s\subseteq M:\{m_1\}\cup s\in P_2 \mbox{ and }\{m_1\}<s\}.
\end{eqnarray*}
Then $\{P_1^1,P_2^1\}$ is a partition of $[M_1]^{<\omega}$. Since
$\ccl(m_1)$ is $\xi_{m_1}$--uniform on $M_1$ and $\xi_{m_1}<\xi$, according
to the induction hypothesis, there exists an infinite subset $L_1$ of $M_1$
such that
$$\ccl(m_1)\cap [L_1]^{<\omega}\subseteq P^1_{i_1}$$
for some $i_1\in \{1,2\}$.

Let $m_2=\min L_1> m_1$ and $M_2=L_1\cap (m_2,+\infty)$. Now set
\begin{eqnarray*}
P^2_1 &=& \{s\subseteq M_2:\{m_2\}\cup s\in P_1\} \mbox{ and }\\
P^2_2 &=& \{s\subseteq M_2:\{m_2\}\cup s\in P_2\}.
\end{eqnarray*}
It is easy to see that $\{P^2_1,P_2^2\}$ is a partition of $M_2$ and that
$\ccl(m_2)\cap [M_2]^{<\omega}$ is $\xi_{m_2}$--uniform on $M_2$
according to Remark 1.2 (ii). Using the induction hypothesis we can find an
infinite subset $L_2$ of $M_2$ and $i_2\in \{1,2\}$ such that
$$\ccl(m_2)\cap [L_2]^{<\omega}\subseteq P^2_{i_2}.$$

Set $m_3=\min L_2$ and $M_3=L_2\cap (m_3,+\infty)$. We proceed inductively
and define a strictly increasing sequence $(m_n)^{\infty}_{n=1}$ in $M$,
two decreasing sequences $(M_n)^{\infty}_{n=1},(L_n)^{\infty}_{n=1}$ in
$[M]$ and a sequence $(i_n)^{\infty}_{n=1}$ in $\{1,2\}$ such that for
every $n\in \nn$ we have
$$m_n=\min L_{n-1},(L_0=M), L_n\subseteq M_n, M_n=L_{n-1}\cap
(m_n,+\infty),\mbox{ and}$$
$$\ccl(m_n)\cap [L_n]^{<\omega}\subseteq P^n_{i_n},$$
where $P^n_{i_n}=\{s\subseteq M_n: \{m_n\}\cup s\in P_{i_n}\}$.

It is clear that there exists an infinite subset $K$ of $\nn$ such that the
subsequence $(i_k)_{k\in K}$ of $(i_n)^{\infty}_n$ is constant; set $i_k=i$
for every $k\in K$, and
$$L=\{m_k:k\in K\}.$$
Then
$$\ccl\cap [L]^{<\omega}\subseteq P_i.$$
Indeed, let $F\in \ccl\cap [L]^{<\omega}$. Then $\min F=m_n$ for some $n\in
K$. Since $F\in \ccl$, we can find $s\in \ccl (m_n)$, such that $m_n<\min s$
and $F=\{m_n\}\cup s$. Also, since $L\cap (m_n,+\infty)\subseteq L_n$ for
every $n\in \nn$, we have that $s\in [L_n]^{<\omega}$. Hence
$$s\in \ccl (m_n)\cap [L_n]^{<\omega}\subseteq P^n_{i_n}.$$
According to the definition of $P^n_{i_n}$ $(n\in \nn)$ we have that $F\in
P_{i_n}$ and, since $n\in K$, that $F\in P_i$.

The proof is complete.

\ \\
\indent
The following Corollary is the precise generalization of the classical
Ramsey theorem.
\begin{coro}
{\rm Let $\xi$ be a countable ordinal, $M\in [\nn]$ and $\ccl$ a
$\xi$--uniform family on $M$. For every $N\in [M]$ and every partition
$\{P_1,\ldots,P_k\}$ of $[N]^{<\omega}$ there exist $L\in [N]$ and $i\in
\{1,\ldots,k\}$ such that:
$$\ccl \cap [L]^{<\omega}\subseteq P_i.$$}
\end{coro}
\pr Let $k=2$. For every $N\in [M]$ the family $\ccl \cap [N]^{<\omega}$ is
$\xi$--uniform on $N$ and \\
$\{P_i\cap [N]^{<\omega}:1\leq i \leq 2\}$ is a
partition of $[N]^{<\omega}$. Hence the proof is immediate by Theorem 2.1.
 
The general case follows by induction on $k$.

\ \\
\indent
In the following corollary we will describe a condition in order for a family
$\ccf$ of finite subsets of $\nn$ to contain a uniform family.
\begin{coro}
{\rm Let $\ccf$ be a family of finite subsets of $\nn$, $M$ an infinite subset
of $\nn$, $\xi$ a countable ordinal and $\ccl$ a $\xi$--uniform family on
$M$. If $\ccl \cap \ccf \cap [L]^{<\omega}\neq \varnothing$ for every $L\in
[M]$, then for every $N\in [M]$ there exists $L\in [N]$ such that:
$$\ccl\cap [L]^{<\omega}\subseteq \ccf.$$}
\end{coro}
\pr Let $N\in [M]$. We set $P_1=\ccf\cap [N]^{<\omega}$ and
$P_2=[N]^{<\omega}\smallsetminus P_1$. According to Theorem 2.2 there
exists $L\in [N]$ such that either $\ccl\cap [L]^{<\omega}\subseteq \ccf$,
as required, or $\ccl\cap [L]^{<\omega}\subseteq
[N]^{<\omega}\smallsetminus \ccf$, which is impossible from our hypothesis.

\ \\
\indent
As we observed in Example 1.12 a uniform family $\ccl$ is not necessary
Sperner. Using Theorem 2.1 it is easy to prove that for every uniform
family $\ccl$ there exists $L\in [\nn]$ such that $\ccl \cap [L]^{<\omega}$
is Sperner.
\begin{coro}
{\rm Let $\xi$ a countable ordinal with $1\leq \xi$, $M\in [\nn]$ and
$\ccl$ a $\xi$--uniform family on $M$. Then there exists $L\in [M]$ such
that $\ccl\cap [L]^{<\omega}$ is Sperner.}
\end{coro}
\pr Let $\ccf =\{ s\in \ccl:$ there is no $t\in \ccl$ such that
$t\subseteq s\}$. It is easy to see that $\ccl \cap \ccf \cap
[L]^{<\omega}\neq \varnothing$ for every $L\in [M]$. Hence according to the
previous corollary there exists $L\in [M]$ such that $\ccl \cap
[L]^{<\omega}\subseteq \ccf$. This gives that $\ccl \cap [L]^{<\omega}$ is
Sperner. 

\ \\
\indent
In the following we will prove that every finite subset of $\nn$ has a
``canonical representation'' with respect to a $\xi$--uniform family on $\nn$,
for every $1\leq \xi <\omega_1$. Using this fact we will prove a 
dichotomy result (Theorem 2.11 below) for hereditary families which is
stronger than Theorem 2.1.
\begin{defi}
{\rm Let $\ccl$ a family of finite subsets of $\nn$ and $A$ a non--empty
finite subset of $\nn$. We will say that $A$ has {\bf canonical
representation} $R_{\ccl}(A)=(s_1,\ldots,s_n,s_{n+1})$, with {\bf type}
$t_{\ccl}(A)=n$, {\bf with respect to ${\bf \ccl}$}, if there exist unique
$n\in \nn$, $s_1,\ldots,s_n\in \ccl$ and $s_{n+1}$ a proper initial segment
of some element of $\ccl$, such that
$A=\displaystyle{\bigcup^{n+1}_{i=1}s_i}$ and $s_1<\ldots<s_n<s_{n+1}$.}
\end{defi}
\begin{prop}
{\rm Let $M$ be an infinite subset of $\nn$, $\xi$ a countable ordinal and
$\ccl$ a $\xi$--uniform family on $M$. Every non--empty finite subset of $M$
has canonical representation with respect to $\ccl$.}
\end{prop}
\pr We will proceed by induction on $\xi$. For $\xi=1$ we have
$\ccl=\{\{m\}:m\in M\}$; if $A=\{m_1,\ldots,m_n\}\in [M]^{<\omega}$, with
$m_1<\ldots<m_n$, then $R_{\ccl}(A)=(\{m_1\},\ldots,\{m_n\})$,
$t_{\ccl}(A)=n$. 

Assume that $1<\xi$ and the assertion holds for every $\zeta<\xi$; and let
$\ccl$ be a $\xi$--uniform family on $M$. Then there exists a sequence
$(\xi_m)_{m\in M}$ of ordinal numbers smaller than $\xi$ such that
$\ccl(m)$ is a $\xi_m$--uniform family on $M\cap (m,+\infty)$ for every
$m\in M$.

Firstly, we will prove that for every $A\in [M]^{<\omega}$, $A\neq
\varnothing$ there exist $n\in \nn$ and $s_1,\ldots,s_n,s_{n+1}\in
[M]^{<\omega}$ such that $A=\displaystyle{\bigcup_{i=1}^{n+1}s_i}$,
$s_1<\ldots<s_n<s_{n+1}$, $s_i\in \ccl$ for every $1\leq i \leq n$ and
$s_{n+1} \prec s_0$ for some $s_0\in \ccl$ (i.e. $s_{n+1}$ as a proper
initial segment of $s_0$).

Let $A\in [M]^{<\omega}$ with $A\neq \varnothing$. If $A\in
\ccl^{*}\smallsetminus \ccl$ then set $n=0$ and $s_1=A$. If $A\in \ccl$
then set $n=1$ and $s_1=A$, $s_2=\varnothing$. So assume that $A\not\in
\ccl^{*}$; then $A=\{m_1\}\cup t^1$ with $t^1\neq \varnothing$ and
$\{m_1\}<t^1$. Since $t^1\in \ccl (m_1)$, according to the induction
hypothesis, $t^1$ has canonical representation
$R_{\ccl(m_1)}(t^1)=(t^1_1,\ldots,t^1_{n_{1}+1})$ with type
$t_{\ccl(m_1)}(t^1)=n_1$ with respect to $\ccl (m_1)$. In this case
$n_1\geq 1$. Indeed, if $n_1=0$, then $A\in \ccl^{*}$, contrary to our
assumption. Set $s_1=\{m_1\}\cup t^1_1$. Obviously, $s_1\in \ccl$,
$s_1\prec A$ and $s_1\neq A$.

We continue analogously setting $A_1=A\smallsetminus s_1$ and treating
$A_1$ in place of $A$ in order to define $s_2$. In detail the argument
goes as follows: if $A_1\in \ccl^{*}\smallsetminus \ccl$ then set $n=1$ and
$s_2=A_1$. If $A_1\in \ccl$ then set $n=2$ and $s_2=A_1$,
$s_3=\varnothing$. Assume that $A_1\not\in \ccl^{*}$; then
$A_1=\{m_2\}\cup t^2$ with $t^2\neq \varnothing$ and $t^2\in \ccl(m_2)$. If
$R_{\ccl(m_2)}(t_2)=(t^2_1,\ldots,t^2_{n_2},t^2_{n_2+1})$ with
$t_{\ccl(m_2)}(t_2)=n_2$, then $n_2\geq 1$. So set $s_2=\{m_2\}\cup t^2_1$
and obviously in this case $s_2\in \ccl$, $s_1\cup s_2\prec A$ and $s_1\cup
s_2\neq A$. Set $A_3=A\smallsetminus s_1\cup s_2$, and continue in the same
way.

Secondly, we prove that for every $A\in [M]^{<\omega}$, $A\neq \varnothing$
the choice of such $n\in \nn$ and sets $s_1,\ldots,s_n,s_{n+1}$ is unique;
so that in fact $t_{\ccl}(A)=n$ and
$R_{\ccl}(A)=(s_1,\ldots,s_n,s_{n+1})$. Indeed, let $A\in [M]^{<\omega}$,
$A\neq \varnothing$ and $t_1,\ldots,t_m,t_{m+1}\in [M]^{<\omega}$ such that
$A=\displaystyle{\bigcup_{i=1}^{m+1}t_i}$, $t_1<\ldots<t_m<t_{m+1}$,
$t_i\in \ccl$ for every $1\leq i\leq m$ and $t_{m+1}\prec t_0$ for some
$t_0\in \ccl$. We will prove, by induction on $m$, that $m=n$ and $t_i=s_i$
for every $1\leq i\leq n+1$.

Let $m=0$. Then $A=t_1$ and there exists $t_0\in \ccl$ such that $t_1\prec
t_0$ and $t_0\neq t_1$. We claim that $n=0$ and consequently $s_1=A=t_1$.
Indeed, if $n\geq 1$ then we have $s_1\prec t_0$, $s_1\neq t_0$ and
$s_1,t_0\in\ccl$, which is impossible since $\ccl$ is thin.

If $m=k+1$ and the assertion holds for $m=k$, then, since $m\geq 1$, we
have, as in the case $m=0$, that $n\geq 1$. Hence, since $t_1
\prec A$, $s_1\prec A$, $t_1,s_1\in \ccl$ and $\ccl$ is thin, we have
that $t_1=s_1$. Set $A_1=A\smallsetminus t_1$; then according to the
induction hypothesis $m=n$ and $t_i=s_i$ for every $1\leq i\leq m+1$.
\begin{coro}
{\rm If $M\in [\nn]$, $\ccl$ is a uniform family on $M$, $M_1\in [M]$, $A$
is a finite subset of $M_1$ and $\ccl_1=\ccl\cap [M_1]^{<\omega}$, then
$t_{\ccl}(A)=t_{\ccl_1}(A)$ and $R_{\ccl}(A)=R_{\ccl_1}(A)$.}
\end{coro}
\pr This holds since the canonical representation of $A$ with respect to
$\ccl$ is unique.

\ \\
\indent
The principal use of the canonical representation of a finite set of $\nn$
in Ramsey theory is contained in the following important Corollary 2.9.
\begin{coro}
{\rm Let $M\in [\nn]$ and $\ccl$ a uniform family on $M$. For every finite,
non empty subset $A$ of $M$ exact one of the following possibilities occurs:\\
either (i) there exists $s\in \ccl$ such that $A\prec s$ and $A\neq s$;\\
or (ii) there exists $s\in \ccl$ such that $s\prec A$.}
\end{coro}
\pr If $A\in [M]^{<\omega}$, $A\neq \varnothing$, then according to
Proposition 2.7, either $t_{\ccl}(A)=0$ (which equivalently gives (i)) or
$t_{\ccl}(A)\geq 1$ (which equivalently gives (ii)).
\begin{coro}
{\rm Let $M\in [\nn]$ and $\ccl$ a uniform family on $M$. If 
$s$ is a proper initial segment of some
element of $\ccl$, then for every $m\in M$ with $s<\{m\}$, the set $s\cup
\{m\}$ is an initial segment of some element of $\ccl$.}
\end{coro}
\pr For every $m\in M$, obviously $\{m\}\in \ccl^{*}$. Let $s\in
\ccl^{*}\smallsetminus \ccl$ with $s\neq \varnothing$ and $m\in M$ with
$s<\{m\}$. Set $A=s\cup\{m\}$. According to Corollary 2.9, either there
exists $s_1\in \ccl$ such that $A\prec s_1$ and $A\neq s_1$ or there
exists $s_2\in \ccl$ such that $s_2\prec A$. In the second case,
we have $s_2=A\in \ccl$, since $\ccl$ is thin. Hence, in both cases 
$A\in \ccl^{*}$.

\ \\
\indent
According to Corollary 2.5 for every uniform family $\ccl$ on $M$ ($M\in
[\nn]$) there exists $L\in [M]$ such that $\ccl\cap [L]^{<\omega}$ is a
Sperner uniform family on $L$. For Sperner uniform families, we have in
fact the following equalities.
\begin{coro}
{\rm Let $M\in [\nn]$ and $\ccl$ a Sperner uniform family on $M$. Then
\begin{enumerate}
\item [(i)] $\ccl^{*}=\ccl_{*}$;
\item [(ii)] $\ccl^{*}\cap [L]^{<\omega}=(\ccl \cap [L]^{<\omega})^{*}$,
and $\ccl_{*}\cap [L]^{<\omega}=(\ccl \cap [L]^{<\omega})_{*}$ for every
$L\in [M]$.
\end{enumerate}}
\end{coro}
\pr \begin{enumerate}
\item [(i)] Obviously, $\ccl^{*}\subseteq \ccl_{*}$. Let $A\in
\ccl_{*}\smallsetminus \ccl$. Then, there exists $s_0\in \ccl$ such that
$A\subseteq s_0$ and $A\neq s_0$. According to Corollary 2.9, either there
exists $s\in \ccl$ such that $A\prec s$ and $A\neq s$ so that $A\in
\ccl^{*}\smallsetminus \ccl$ ensues or there exists $s\in \ccl$ such that
$s\prec A$, an imposibility, since $s\subseteq s_0$ and $\ccl$ is Sperner.
\item [(ii)] Obviously, $(\ccl\cap [L]^{<\omega})^{*}\subseteq \ccl^{*}\cap
[L]^{<\omega}$, and according to (i) $(\ccl\cap [L]^{<\omega})_{*}=
(\ccl\cap [L]^{<\omega})^{*}$, for every $L\in [M]$. Let $L\in [M]$ and $A\in
\ccl^{*}\cap [L]^{<\omega}\smallsetminus \ccl$. According to Corollary
2.10 there exists $s_0\in \ccl$ such that $A\prec s_0$ and
$s_0\subseteq L$. Hence, $A\in (\ccl\cap [L]^{<\omega})^{*}$. This
establishes the required equalities. 
\end{enumerate}

Using Corollary 2.9 (to the canonical representation of finite subsets of
$\nn$) and the general Ramsey theorem (Theorem 2.2) we now prove a stronger
dichotomy result for hereditary families.
\begin{theo}
Let $M\in [\nn]$, $\ccf$ a hereditary family of finite subsets of $M$ and
$\ccl$ a uniform family on $M$. Then for every $M_1\in [M]$ there exists
$L\in [M_1]$ such that 
either $\ccl_{*}\cap [L]^{<\omega}\subseteq \ccf$, 
or $\ccf \cap [L]^{<\omega}\subseteq \ccl^{*}\smallsetminus \ccl$.
\end{theo}
\pr According to Corollary 2.5, there exists $N\in [M_1]$ such that $\ccl
\cap [N]^{<\omega}$ is Sperner. Hence, using Corollary 2.3 of the
general Ramsey theorem (Theorem 2.2) we can find $L\in [N]$ such that 
$$\mbox{either }\ \ccl \cap [L]^{<\omega} \subseteq \ccf, \ \mbox{ or }\
\ccl \cap [L]^{<\omega}\subseteq [\nn]^{<\omega}\smallsetminus \ccf.$$

If $\ccl \cap [L]^{<\omega} \subseteq \ccf$ then $(\ccl \cap
[L]^{<\omega})_{*} \subseteq \ccf$, since $\ccf$ is hereditary. Hence, 
$\ccl_{*} \cap [L]^{<\omega} \subseteq \ccf$, according to Corollary 2.11.

If $\ccl \cap [L]^{<\omega} \subseteq [\nn]^{<\omega}\smallsetminus \ccf$,
we will prove that $\ccf \cap [L]^{<\omega}\subseteq \ccl^{*}\smallsetminus
\ccl$. Indeed, let $A\in \ccf \cap [L]^{<\omega}$. Since $\ccl \cap
[L]^{<\omega}$ is uniform on $L$, according to Corollary 2.9, either there exists
$s\in \ccl \cap [L]^{<\omega}$ such that $A\prec s$ and $A\neq s$, which
gives that $A\in \ccl^{*}\smallsetminus \ccl$, as required, or there exists
$s\in \ccl \cap [L]^{<\omega}$ such that $s\prec A$. But this case is
impossible, since then $s\in \ccf \cap \ccl \cap [L]^{<\omega}$. This
completes the proof.
\begin{coro}
{\rm Let $M\in [\nn]$, $\ccl_1$ a $\xi_1$--uniform on $M$ and $\ccl_2$ a
$\xi_2$--uniform on $M$ with $\xi_1<\xi_2<\omega_1$. Then there exists
$L\in [M]$ such that 
$$\ccl_1\cap [L]^{<\omega}\subseteq (\ccl_2)^{*}\smallsetminus \ccl_2.$$}
\end{coro}
\pr According to Theorem 2.12 there exists $L\in [M]$ such that, 

either $(\ccl_2)_{*}\cap [L]^{<\omega}\subseteq (\ccl_1)_{*}$, 
or $(\ccl_1)_{*}\cap [L]^{<\omega}\subseteq (\ccl_2)_{*}\smallsetminus
\ccl_2.$

The first alternative, is impossible; in fact, if (i) holds, then\\
$s_L((\ccl_2)_{*}\cap [L]^{<\omega})\leq s_L ((\ccl_1)_{*}\cap
[L]^{<\omega})$. On the other hand, using Remark 1.11 (iii) and Theorem 1.18, we
have 
$$s_L((\ccl_2)_{*}\cap [L]^{<\omega})= s_L ((\ccl_2)_{*})=\xi_2+1$$
and
$$s_L((\ccl_1)_{*}\cap [L]^{<\omega})= s_L ((\ccl_1)_{*})=\xi_1+1<\xi_2+1.$$
A contradiction. Thus (ii) holds, as required.

\ \\
\indent
In the following theorem (Theorem 2.15) we will describe, with the help of the strong
Cantor Bendixson index, sufficient conditions in order a family of finite
subsets of $\nn$ to satisfy
exactly one of the conditions given in the dichotomy of Theorem 2.12.

Since we will restrict to the hereditary and closed families firstly we
will give a characterization of them.
\begin{prop}
{\rm Let $\ccf$ be a non empty, hereditary family of finite subsets of $\nn$.
The following are equivalent:
\begin{enumerate}
\item [(i)] $\ccf$ is closed.
\item [(ii)] There does not exist an infinite sequence
$(s_i)^{\infty}_{i=1}$ of elements of $\ccf$ with $s_1\prec s_2 \prec
\ldots$. 
\item [(iii)] There does not exist $M\in [\nn]$ such that
$[M]^{<\omega}\subseteq \ccf$.
\end{enumerate}}
\end{prop}
\pr (i) $\Rightarrow$ (ii) If $s_i=(n_1,\ldots,n_{k_i})\in \ccf$ with
$n_1<\ldots<n_{k_i}$, for every $i\in \nn$ and $(k_i)^{\infty}_{i=1}$ is an
increasing sequence of natural numbers, then $(s_i)^{\infty}_{i=1}$
converges pointwise to the infinite subset $s=(n_1,n_2,\ldots)$ of $\nn$
which does not belong to $\ccf$.

\smallskip
\indent
(ii) $\Rightarrow$ (i) Let $(t_n)^{\infty}_{n=1}$ a sequence of elements of
$\ccf$, converging pointwise to some subset $t$ of $\nn$.
If $t$ is finite, then there exists $n_0\in \nn$ such that
$t\prec t_{n_0}$, hence $t\in \ccf$, as required.

Let $t$ is infinite. Set $t=(n_1,n_2,\ldots)$ with
$n_1<n_2<\ldots$ and $s_i=(n_1,\ldots,n_i)$ for every $i\in \nn$. For
every $i\in \nn$ the
sequence $(t_n\cap [0,n_i])^{\infty}_{n=1}$ converges pointwise to $s_i$. 
According to the previous case, we have $s_i\in \ccf$
for every $i\in \nn$. A contradiction to the condition (ii).

\smallskip
\indent
(iii) $\Rightarrow$ (ii) Let $(s_i)^{\infty}_{i=1}\subseteq \ccf$ with
$s_1\prec s_2\prec\ldots$. Set $\displaystyle{M=\bigcup^{\infty}_{i=1}s_i}$.
If $t$ is an arbitrary subset of $M$ then $t\subseteq s_i$ for some $i\in
\nn$, hence $t\in \ccf$. So $[M]^{<\omega}\subseteq \ccf$, a contradiction. 

\smallskip
\indent
(ii) $\Rightarrow$ (iii) Let $M=(m_1,m_2,\ldots)\subseteq \nn$ with
$m_1<m_2<\ldots$. If $[M]^{<\omega}\subseteq \ccf$  then 
$s_i=(m_1,\ldots,m_i)\in \ccf$ for every $i\in \nn$. Hence, the condition
(ii) does not hold.

\ \\
\indent
After the previous proposition we can give a dichotomy result rather closed to the
infinite Ramsey theorem (c.f. Nash--Williams \cite{11}, Galvin--Prikry
\cite{8}, Silver \cite{16}), and in many (especially Banach space-)
applications it can be used in its place.
\begin{coro}
{\rm Let $\ccf$ be a hereditary family of finite subsets of $\nn$. For every
$M\in [\nn]$ there exists $L\in [M]$ such that either
$[L]^{<\omega}\subseteq \ccf$ or $[L]^{<\omega}\subseteq
([\nn]^{<\omega}\smallsetminus \ccf)_{*}$.}
\end{coro}
\pr According to
Proposition 2.14, if $\ccf\cap [M]^{<\omega}$ is not closed then there
exists $L\in [M]$ such that $[L]^{<\omega}\subseteq\ccf$ and 
if $\ccf\cap [M]^{<\omega}$ is closed then, there is $L\in [M]$ such that
$[L]^{<\omega}\subseteq ([\nn]^{<\omega}\smallsetminus \ccf)_{*}$.

\begin{theo}
Let $\ccf$ be a pointwise closed and hereditary family of finite subsets of
$\nn$, $M\in [\nn]$, $\xi$ a countable ordinal and $\ccl$ a $\xi$--uniform
family on $M$.
\begin{enumerate}
\item [(i)] If $\xi+1<s_M(\ccf)$, then there exists $L\in [M]$ such that 
$$\ccl_{*}\cap [L]^{<\omega}\subseteq \ccf;\mbox{ and},$$
\item [(ii)] If $s_M(\ccf)<\xi+1$, then there exists $L\in [M]$ such that
$$\ccf\cap [L]^{<\omega}\subseteq \ccl^{*}\smallsetminus \ccl.$$
\end{enumerate}
\end{theo}
\pr Using Theorem 2.12 for the family $\ccf\cap [M]^{<\omega}$, (at least)
one of the following possibilities occurs: either there exists $L\in [M]$
such that $\ccl_{*}\cap [L]^{<\omega}\subseteq \ccf$ or there exists $L\in
[M]$ such that $\ccf\cap [L]^{<\omega}\subseteq \ccl^{*}\smallsetminus \ccl$.
\begin{enumerate}
\item [(i):] If $\xi+1<s_M(\ccf)$, then the second case cannot occur, since
then, we would have\\
$\ccf\cap [L]^{<\omega}\subseteq \ccl^{*}\cap [L]^{<\omega}$, and
consequently, according to Theorem 1.18, 
$$\xi+1=s_L(\ccl^{*})=s_L(\ccl^{*}\cap [L]^{<\omega})\geq s_L(\ccf \cap
[L]^{<\omega})=s_L(\ccf),$$ 
a contradiction. 
\item [(ii):] If $s_M(\ccf)<\xi+1$ then the first case can not occur, since
then
$$\xi+1=s_L(\ccl^{*})\leq s_L(\ccf),$$
a contradiction to our hypothesis.
\end{enumerate}
\begin{rem2}
{\rm It should be noted that in the limiting case $s_M(\ccf)=\xi+1$ of Theorem
2.16 both alternatives may materialize. Indeed, we have the following two
simple examples:}
\end{rem2}

\smallskip
\noindent
{\bf Example 1.} Let
$$\ccl=\{s\in [\nn]^{<\omega}:|s|=2\min s+1\}\ \mbox{ and}$$
$$\ccr=\{s\in [\nn]^{<\omega}:|s|=\min s\},$$
where $|s|$ denotes the cardinality or $s$.

It is easy to see that $\ccl$ and $\ccr$ are $\omega$--uniform on $\nn$.

The family $\ccf_1=\ccr_{*}$ is hereditary, closed (Lemma 1.17) and
$s_{\nn}(\ccf)=\omega+1$, according to Theorem 1.18. Since $\ccl\cap
\ccf_1=\varnothing$ and $\ccl\cap [L]^{<\omega}\neq \varnothing$ for every
$L\in [\nn]$ the first alternative of Theorem 2.16 does not occur. Hence
there exists $L\in [\nn]$ such that
$$\ccf_1\cap [L]^{<\omega}\subseteq \ccl^{*}\smallsetminus \ccl.$$

\smallskip
\noindent
{\bf Example 2.} On the other hand (refering the notation of Example 1) for
the hereditary and closed family $\ccf_2=\ccl_{*}$ with
$s_{\nn}(\ccf_2)=\omega+1$ (Theorem 1.18) and the $\omega$--uniform family
$\ccr$ on $\nn$ we have that
$$\ccr^{*}\subseteq \ccf_2$$
hence the first alternative of Theorem 2.16 occurs and the second does not
occur since for every $L\in [\nn]$ there exists $s\in \ccf_2\cap
[L]^{<\omega}$ such that $s\not\in \ccr^{*}$ (take $s\in [L]^{<\omega}$
with $\min s+1\leq |s|\leq 2\min s+1$).

\ \\
\noindent
{\Large\bf Recapitulation of the main results}

\smallskip
\indent
Let $\ccf$ a hereditary family of finite subsets of $\nn$. We have the
following two cases:

\smallskip
\noindent
{\bf 1st case.} The family $\ccf$ is not closed. Then according to
Proposition 2.14 there exists $L\in [\nn]$ such that
$[L]^{<\omega}\subseteq \ccf$.

\smallskip
\noindent
{\bf 2nd case.} The family $\ccf$ is closed. Then there exists $L\in
[\nn]$ such that $[L]^{<\omega}\subseteq ([\nn]^{<\omega}
\smallsetminus\ccf)_{*}$, (Corollary 2.15). Moreover, for a
given infinite subset $M$ of $\nn$ and a system of uniform families 
$(\cca_{\zeta})_{\zeta<\omega_1}$ on $M$, setting
$$\xi =\sup\{s_L(\ccf):L\in [M]\}$$
the following obtain:
\begin{enumerate}
\item [(i)] For every ordinal $\zeta$ with $\zeta+1<\xi$ there exists $L\in
[M]$ such that:
$$(\cca_{\zeta})_{*}\cap [L]^{<\omega}\subseteq \ccf;$$
(Theorem 2.16).
\item [(ii)] For every ordinal $\zeta$ with $\xi<\zeta+1$ and for every
$M_1\in [M]$ there exists $L\in [M_1]$ such that:
$$\ccf\cap [L]^{<\omega}\subseteq (\cca_{\zeta})^{*}\smallsetminus
\cca_{\zeta}$$ 
which gives that
$$\cca_{\zeta}\cap [L]^{<\omega}\subseteq [\nn]^{<\omega}\smallsetminus
\ccf$$ 
(Theorem 2.16).
\item [(iii)] If $\xi=\zeta+1$ then there exists $M_1\in [M]$ such that
$s_{M_1}(\ccf)=\xi$. According to Theorem 2.12 there exists $L\in
[M_1]$ such that

either $(\cca_{\zeta})_{*}\cap [L]^{<\omega}\subseteq \ccf$, 
or $\ccf \cap [L]^{<\omega}\subseteq (\cca_{\zeta})^{*}\smallsetminus
\cca_{\zeta}.$ 
\end{enumerate}

According to Remark 2.17 both alternatives may materialize.

\section{Some remarks and applications}

The results of the previous section constitute a far reaching and powerful
generalization of the classical Ramsey theorem, a generalization that is
stated in terms of a countable ordinal index $\xi$ (in place of a natural
number as in the classical case); these ordinal index dichotomies are in
turn analogous to the Galvin--Prikry (\cite{8}) infinitary form of the
Ramsey dichotomy (stated for all infinite subsets of $\nn$, partitioned by
an analytic partition). Our results, then, on the one hand generalize the
classical Ramsey theorem and on the other hand they have the Galvin--Prikry
infinitary theorem as the limiting $(\omega_1-)$ case.

It is to be expected that such general combinatorial principles will have
wide applications in many instances, where either classical Ramsey theory,
or the infinitary Galvin--Prikry theorem has been successfully applied.
Some applications of the dichotomy results established in this paper, and
their relation to existing applications of similar combinatorial techniques
involving mostly generalized Schreier families $F_{\xi},\xi<\omega_1$ (as
in \cite{2}, \cite{6}, \cite7{}, \cite{10}), will appear in a separate
publication. 

Here we will limit ourselves to exhibit the way in which our techniques can
be applied to provide simple derivations of the combinatorial basis in the
theory of Banach spaces of two recent results, one by
Argyros--Mercourakis--Tsarparlias and the other by Judd.

Thus in Proposition 3.1 we indicate the close connection that exists
between the Schreier family $F_{\xi}$ and an $\omega^{\xi}$--uniform
family, particularly the family $\ccb_{\xi}$ (Definition 1.5). Then using the
results of Section 2 we reprove a dichotomy result of Judd \cite{9},
obtaining in fact a more general expression; and additionally a
combinatorial result of Argyros--Mercourakis--Tsarpalias \cite{2} which was
the basis for establishing a general form of an $\ell^1$--dichotomy,
initially proved in a special form by Rosenthal \cite{15}.
\begin{prop}
{\rm Let $\xi$ be a countable ordinal, $M$ an infinite subset of $\nn$ and
$\ccl$ an $\omega^{\xi}$--uniform family on $M$. Then there exists $L\in
[M]$ such that $\ccf_{\xi}(L)\subseteq \ccl^{*}$.}
\end{prop}
\pr Let $\ccl^1=\{\{m\}\cup s:m\in M, s\in \ccl \mbox{ and }\{m\}<s\}$.
According to Remark 1.2 (iii) $\ccl^1$ is $\omega^{\xi}+1$--uniform on
$M$. 
Using Theorem 2.12, either there exists $N\in [M]$ such that \\
$(\ccl^1)_{*}\cap [N]^{<\omega}\subseteq \ccf_{\xi}$, which is impossible, since
$s_N[(\ccl^1)_{*}]=\omega^{\xi}+2$ (according to Theorem 1.18) and
$s_N(\ccf_{\xi})=\omega^{\xi}+1$ (according to \cite{2}); or there exists
$N\in [M]$ such that $\ccf_{\xi}\cap [N]^{<\omega}\subseteq
(\ccl^1)^{*}\smallsetminus \ccl^1$. In this case, set
$L=(n_i)^{\infty}_{i=3}$ if $N=(n_i)^{\infty}_{i=1}$. 
Then it is easy to see that $\ccf_{\xi}(L)\subseteq \ccl^{*}$, using the
fact that if $(k_1,\ldots,k_p)\in \ccf_{\xi}$ then
$(k_1+1,k_1+2,k_2+2,\ldots,k_p+2)\in \ccf_{\xi}$ for every $\xi<\omega_1$.
\begin{coro}
{\rm For every $\xi <\omega_1$, $M\in [\nn]$, there exists $L\in [M]$ such
that 
$$\ccf_{\xi}(L)\subseteq (\ccb_{\xi})^{*}\subseteq \ccf_{\xi}$$}
\end{coro}
\pr It is immediate after Theorem 1.6 and Proposition 3.1

\ \\
\indent
After Proposition 3.1 we will give a Corollary of the general Ramsey
Theorem 2.2 which can be used for the families $\ccf_{\xi},\xi<\omega_1$.
\begin{prop}
{\rm Let $\ccf$ be a hereditary family of finite subsets of $\nn$, $M$ an
infinite subset of $\nn$ and $\xi$ a countable ordinal number. If
$\ccf\cap \ccb_{\xi}\cap [L]^{<\omega}\neq \varnothing$ for every $L\in
[M]$ then there exists $L\in [M]$ such that $\ccf_{\xi}(L)\subseteq \ccf$.}
\end{prop}
\pr According to Corollary 2.4 and since $\ccb_{\xi}\cap [M]^{<\omega}$ is
$\omega^{\xi}$--uniform family on $M$ (Corollary 1.8) there exists $N\in
[M]$ such that $\ccb_{\xi}\cap [N]^{<\omega}\subseteq \ccf$. From the
previous proposition there exist $L\in [N]$ such that
$\ccf_{\xi}(L)\subseteq (\ccb_{\xi}\cap [N]^{<\omega})_{*}$. 
Since $\ccf$ is hereditary and $\ccb_{\xi}\cap [N]^{<\omega}\subseteq
\ccf$ we have that $(\ccb_{\xi}\cap [N]^{<\omega})_{*}\subseteq \ccf$. 
Hence, $\ccf_{\xi}(L)\subseteq \ccf$, as required.

\ \\
\indent
R. Judd in \cite{9} had provided, using Schreier games, that for every
hereditary family $\ccf$ of finite subsets of $\nn$, $\xi<\omega_1$ and
$M\in [\nn]$, either there exists $L\in [M]$ such that
$\ccf_{\xi}(L)\subseteq \ccf$ or there exists $L\in [M]$ and $N\in [\nn]$
such that $\ccf\cap [N]^{<\omega}(L)\subseteq \ccf_{\xi}.$

We will prove a stronger version of this result using our results of
Section 2.
\begin{theo}
For every hereditary family $\ccf$ of finite subsets of $\nn$, every
countable ordinal $\xi$ and $M\in [\nn]$ there exists $L\in [M]$ such that
either $\ccf_{\xi}(L)\subseteq \ccf$ or $\ccf\cap [L]^{<\omega}\subseteq
\ccf_{\xi}$. 
\end{theo}
\pr According to Theorem 2.12 there exist $N\in [M]$ such that: 

either $(\ccb_{\xi})_{*}\cap [N]^{<\omega}\subseteq \ccf$, 
or $\ccf \cap [N]^{<\omega}\subseteq (\ccb_{\xi})^{*}$.

\indent
Using Proposition 3.1 there exists $L\in [N]$ such that
$$\ccf_{\xi}(L)\subseteq (\ccb_{\xi})^{*}\cap [L]^{<\omega}\subseteq
(\ccb_{\xi})_{*} \cap [N]^{<\omega}.$$
Hence, there exists $L\in [M]$ such that 
either $\ccf_{\xi}(L)\subseteq \ccf$, 
or $\ccf \cap [L]^{<\omega}\subseteq (\ccb_{\xi})^{*}\subseteq \ccf_{\xi}.$

\ \\
\indent
As a corollary of Theorem 3.3 we have the following result of Argyros,
Mercourakis and Tsarpalias in \cite{2}. An analogous proof for this result
was given in \cite{9}.
\begin{theo}
Let $\ccf$ be a hereditary and closed family of finite subsets of $\nn$.
If there exists $M\in [\nn]$ such that $s_M(\ccf)\geq\omega^{\xi}$, then there
exists $L\in [M]$ such that $\ccf_{\xi}(L)\subseteq \ccf$.
\end{theo}
\pr If $s_M[\ccf]>\omega^{\xi}+1$, then according to Theorem 2.16 (i), there exists
$N\in [M]$ such that\\
 $(\ccb_{\xi})_{*}\cap [N]^{<\omega}\subseteq \ccf$.
Also, according to Proposition 3.1 there exists $L\in [N]$ such that
$$\ccf_{\xi}(L)\subseteq (\ccb_{\xi})^{*}\cap [L]^{<\omega}$$
Hence, $\ccf_{\xi}(L)\subseteq \ccf$.

Now, if $s_M[\ccf]=\omega^{\xi}+1$ then set $\overline{\ccf}=\{\{m\}\cup
s:s\in \ccf,m \in M\mbox{ and }\{m\}<s\}$. It is easy to see that
$s_M[\overline{\ccf}]>\omega^{\xi}+1$. If we apply the previous case to
$\overline{\ccf}$ we can find $(n_i)^{\infty}_{i=1}=N\in [M]$ such that
$\ccf_{\xi}(N)\subseteq \overline{\ccf}$ and setting
$L=(n_i)^{\infty}_{i=3}$ we have that $\ccf_{\xi}(L)\subseteq \ccf$ as
required. 

\ \\
\noindent
{\large\bf Acknowledgement}

\smallskip
\indent
I would like to thank Professor S. Negrepontis and A. Katavolos for their
helpful comments and support during the preparation of this paper.

{\sc Department of Mathematics, Panepistimiopolis, 15784, Athens, Greece.}

{\it E--mail address:} vgeorgil@atlas.uoa.gr
\end{document}